\documentclass[11pt]{article}
\usepackage{graphicx}
\usepackage{mathrsfs,amsmath,amssymb}
\usepackage{commath}
\usepackage[nottoc]{tocbibind}
\usepackage{fullpage}

\usepackage{array}
\usepackage{diagbox} 
\usepackage{caption}
\numberwithin{equation}{section}  
\usepackage{epsfig,amsbsy,amsthm}
\usepackage{float,epsfig}
\usepackage{fixmath,textcomp}
\usepackage{graphics,caption,subcaption}
\usepackage{pgfpages,authblk}
\usepackage{appendix}
\usepackage[numbers,sort&compress]{natbib}
\usepackage{tabularx,ragged2e,caption}
\usepackage{hyperref} 
\usepackage{graphicx}
\usepackage{pst-all}
\allowdisplaybreaks 

\begin{document}
\newtheorem{theorem}{\bf Theorem}[section]
\newtheorem{proposition}[theorem]{\bf Proposition}
\newtheorem{definition}{\bf Definition}[section]
\newtheorem{corollary}[theorem]{\bf Corollary}
\newtheorem{example}[theorem]{\bf Example}
\theoremstyle{remark} 
\newtheorem{remark}[theorem]{\bf Remark}
\newtheorem{lemma}[theorem]{\bf Lemma}
\newtheorem{algorithm}[theorem]{Algorithm}
\newtheorem{observation}[theorem]{Observation}
\newcommand{\von}{\vskip 1ex}
\newcommand{\vone}{\vskip 2ex}
\newcommand{\vtwo}{\vskip 4ex}
\newcommand{\ds}{\displaystyle}
\newcommand*{\rom}[1]{\expandafter\@slowromancap\romannumeral #1@}

\def \noin{\noindent}

\title{{\bf A Biologically Motivated Finite Difference Approach for Simulating Singularly Perturbed Vertical Motion in Human Gait}}

\author[1]{Shubhangini Gupta}
\author[2]{Sourav Banerjee}
\author[3,$\ast$]{Tamal Pramanick}
\affil[1,2,3]{\em \scriptsize {Department of Mathematics, National Institute of Technology Calicut, Calicut, Kerala, India.}}
\affil[1]{\em \scriptsize {email id: shubhangini\_p240187ma@nitc.ac.in}}
\affil[2]{\em \scriptsize {email id: souravsbbanerjee@gmail.com}}
\affil[3]{\em \scriptsize {email id: tamal@nitc.ac.in, $^\ast$Corresponding author}}

\date{}
\maketitle

{\small \textbf{Abstract.} In this study, we present a simulation-based numerical method for solving a class of singularly perturbed second-order differential equations that come from a simplified biologically motivated model of human gait. Important physical factors such as gravity, damping, and leg stiffness are included in the model, which also  depicts the vertical motion of the center of mass of the body during walking or running. Most of the time, standard numerical methods are ineffective in resolving boundary layer behavior that occurs due to the small perturbation parameter in the governing equation. We use a domain decomposition technique to divide the problem domain into inner and outer regions to tackle this difficulty. The boundary layer resolves the steep gradients. We applied a time-rescaling transformation to the inner region. Each subdomain is discretized, and the resulting tridiagonal systems are efficiently solved using the Thomas algorithm within the mixed finite difference framework. A detailed convergence analysis demonstrates second-order accuracy in space.          The numerical results validate the proposed scheme’s accuracy, stability, and efficiency through experiments based on modified human gait models. The framework serves as a fundamental tool for biomechanical simulation. The modeling is a foundation for future research, incorporating nonlinearities, time delays, and real-world scenarios data on how people walk.


}

\section{Introduction}
Mathematical modeling of human gait has emerged as a critical area of research due to its wide-ranging applications in clinical gait analysis, rehabilitation engineering, prosthetic limb development, humanoid robotics, and sports science. The most important aspect of this modeling is the vertical displacement of the body's center of mass (COM), which is crucial for the energy efficiency and mechanical stability of locomotion. Simple mechanical models, like spring-mass systems, inverted pendulum models, and variations on these, have been used to successfully model key aspects of gait mechanics, such as impact forces, energy transfer, leg compliance, and ground response forces, see e.g., \cite{Mochon1980, leveque2007}.

Often using minor physical quantities as stiffness ratios, mass inertia, or damping coefficients, biomechanical models produce governing equations (cf. \cite{duressa2015, holmes2012}). These small parameters cause singular perturbations that occurs due to narrow areas near domain boundaries where the solution experiences fast variation to develop, which is named as boundary layers. Standard numerical approaches, which frequently fail to resolve the steep gradients precisely and may suffer from instability or loss of precision, find difficulty in these layers.  When one models the motion of the center of mass (COM) while walking or running, where the dynamics require a complicated interaction between position, velocity and acceleration. One can readily model such biological processes as periodic dynamic systems.  Often found near the ends of the domain, the resulting singularly perturbed differential or differential-difference equations are notable for their abrupt boundary layer behavior.  Precisely capturing these characteristics calls for specific numerical techniques since ordinary discretization approaches often fail to preserve stability or provide sufficient accuracy in the presence of these layers \cite{Roos2008, kadalbajoo2008, miller1996}.

In this study, we introduce a novel mixed finite difference approach to address the solution of a singularly perturbed differential-difference model that characterizes the intricacies of human gait. Our formulation is inspired by reduced-order nonlinear oscillatory models in biomechanics (see \cite{perc2005, Schmalz2021}) and integrates structural insights from research on singularly perturbed delay differential equations (cf. \cite{Sirisha2018, rai2011,rai2012, Nageshwar2013}). Our method diverges from shift operator-based or solely fitted schemes by using a terminal boundary point to partition the domain into separate inner and outer subregion.  This enables us to convert the initial problem into an asymptotically comparable singular perturbation formulation. The use of time rescaling in the boundary layer area is a key new feature of our framework.  This higher time resolution is especially useful for apps like tracking athletic motion, because it lets us get a clear picture of how quickly the solution changes without having to refine the mesh too much.  In both subdomains, central difference approximations are used, and the Thomas algorithm is used to quickly solve the discrete systems that are left over, see e.g. \cite{song2016, stergiou2004, Mochon1980}.

The vertical motion of the human body's center of mass during locomotion can be described by a second-order ordinary differential equation that incorporates key biomechanical properties such as leg stiffness and damping. Taking into account the small perturbation parameter $\varepsilon$, the governing equation is given by:
\[
z_{tt} + \varepsilon u(t) z_t + v(t) z = -g,
\]
where $u(t) = \varepsilon \mu(t)$ represents a damping function and $v(t)$ denotes the stiffness. The boundary conditions are $z(0) = f_1$, $z(t_f) = f_2$. In the presence of the small parameter $\varepsilon$, the solution exhibits a boundary layer near the domain endpoints.

We use the suggested technique over several terminal point selections and validate its performance on several model instances, including those with perturbation- delay and advancement parameters (cf. \cite{kadalbajoo2008}). 
This phenomena guarantees the computational robustness \cite{smith1985}. We show that the method is well-suited for stiff and layer-dominated regimes generally seen in biomechanical systems since it achieves second-order precision and uniform convergence with respect to the perturbation parameter \( \varepsilon \).

To the best of authors knowledge, this is the first work to introduce a hybrid finite difference framework that combines domain decomposition, boundary layer scaling, and structural adaptation for solving singularly perturbed models in human gait. Numerical experiments are presented to validate the theoretical convergence behavior and demonstrate the method’s versatility and effectiveness. The work is arranged as follows. The mathematical model for the human gait is presented in the next section, namely Section \ref{eq:a0}, together with a time-rescaling technique and domain decomposition analysis. The discretization using the mixed finite difference approach and the algorithmic implementation is covered in Section \ref{eq:a1}, and Section \ref{eq:a2} offers a theoretical convergence study. Numerical experiments verifying the accuracy and homogeneity of convergence of the technique are found in Section \ref{eq:a3}, and section \ref{eq:a5} ends the work with a summary and recommendations for next work including extensions to nonlinear models with changing coefficients and control inputs.

\section{Mathematical Modeling and Problem Formulation}\label{eq:a0}
The vertical motion of the human body's center of mass during locomotion can be modeled using a second-order singularly perturbed differential equation incorporating biomechanical effects such as damping and leg stiffness. The governing equation is given by:
\begin{equation}
    z_{tt} + \varepsilon u(t) z_t + v(t) z = -g, \label{eq:gov}
\end{equation}
where
   $\varepsilon$ is a small positive perturbation parameter ($0 < \varepsilon \ll 1$),
     $u(t) = \varepsilon \mu(t)$ is the damping function,
     $v(t)$ represents leg stiffness,
     $g$ is the gravitational constant,
   $z(t)$ is the vertical displacement of the COM.
The boundary conditions are:
\begin{equation*}
    z(0) = f_1, \quad z(t_f) = f_2. \label{eq:bc}
\end{equation*}

\subsection{Time Rescaling in the Boundary Layer}
To resolve the steep gradients near the boundary at $t=0$, we introduce a stretched variable:
\begin{equation*}
    T = \varepsilon t \quad \Rightarrow \quad t = \frac{T}{\varepsilon}.
\end{equation*}
Let $Z(T) = z\left(\frac{T}{\varepsilon}\right)$, then the derivatives transform as:
\begin{equation*}
    z_t = \varepsilon Z_T, \qquad z_{tt} = \varepsilon^2 Z_{TT}.
\end{equation*}
Substituting into equation~\eqref{eq:gov}, we obtain the transformed equation in the inner region:
\begin{equation}
    \varepsilon^2 Z_{TT}(T) + \varepsilon^2 U(T) Z_T(T) + V(T) Z(T) = -g, \label{eq:inner}
\end{equation}
where $U(T) = \mu(T/\varepsilon)$ and $V(T) = v(T/\varepsilon)$, with modified boundary conditions:
\begin{equation*}
    Z(0) = f_1, \quad Z(\varepsilon t_f) = f_2.
\end{equation*}
To accurately solve the singularly perturbed second-order differential equation representing human gait dynamics, we employ a mixed finite difference method based on domain decomposition and time rescaling. This method divides the domain into two regions, inner and outer regions. On the inner region a steep gradients (boundary layers) are expected to occur and a stretched variable is used to resolve this. On the other hand, the outer region, where the solution varies smoothly, standard discretization is sufficient.

\section{Domain Decomposition using Mixed Finite Difference Scheme}\label{eq:a1}
We divide the computational domain $[0, t_f]$ into two subdomains to address boundary layer behavior:
\begin{itemize}
    \item \textbf{Inner region}: $T \in [0, T_p]$, where $T_p = \varepsilon t_p$ captures the boundary layer, and
    \item \textbf{Outer region}: $t \in [t_p, t_f]$, where the solution varies smoothly.
\end{itemize}

\subsection{Discretization in the Inner Region}
Let $T_i = i k$ be a uniform discretization of the interval $[0, T_p]$, where $k$ is the mesh step size. Denoting $Z(T_i) = Z_i$, the central difference approximations are:
\begin{align*}
    Z_{TT}(T_i) &\approx \frac{Z_{i+1} - 2Z_i + Z_{i-1}}{k^2}, \\
    Z_T(T_i) &\approx \frac{Z_{i+1} - Z_{i-1}}{2k}.
\end{align*}
Substituting these into equation \eqref{eq:inner}, we obtain the finite difference scheme:
\begin{equation*}
    \varepsilon^2 \left( \frac{Z_{i+1} - 2Z_i + Z_{i-1}}{k^2} \right) + \varepsilon^2 U_i \left( \frac{Z_{i+1} - Z_{i-1}}{2k} \right) + V_i Z_i = -g,
\end{equation*}
which can be rearranged as:
\begin{equation*}
    C_i Z_{i-1} + A_i Z_i + B_i Z_{i+1} = -g,
\end{equation*}
where:
\begin{align*}
    C_i &= \varepsilon^2 \left( \frac{1}{k^2} + \frac{U_i}{2k} \right), \\
    A_i &= V_i - \frac{2\varepsilon^2}{k^2}, \\
    B_i &= \varepsilon^2 \left( \frac{1}{k^2} - \frac{U_i}{2k} \right).
\end{align*}

\subsection{Discretization in the Outer Region}
Similarly, for the outer region $t \in [t_p, t_f]$ with uniform step size $k$ and grid points $t_i = t_p + i k$, denoting $z(t_i) = z_i$, we use:
\begin{align*}
    z_{tt}(t_i) &\approx \frac{z_{i+1} - 2z_i + z_{i-1}}{k^2}, \\
    z_t(t_i) &\approx \frac{z_{i+1} - z_{i-1}}{2k}.
\end{align*}
Substituting into the original equation~\eqref{eq:gov} gives:
\begin{equation*}
    \frac{z_{i+1} - 2z_i + z_{i-1}}{k^2} + \varepsilon u_i \cdot \frac{z_{i+1} - z_{i-1}}{2k} + v_i z_i = -g,
\end{equation*}
which we write as:
\begin{equation*}
    \widetilde{C}_i z_{i-1} + \widetilde{A}_i z_i + \widetilde{B}_i z_{i+1} = -g,
\end{equation*}
where:
\begin{align*}
    \widetilde{C}_i &= \frac{1}{k^2} - \frac{\varepsilon u_i}{2k}, \\
    \widetilde{A}_i &= v_i - \frac{2}{k^2}, \\
    \widetilde{B}_i &= \frac{1}{k^2} + \frac{\varepsilon u_i}{2k}.
\end{align*}

These two systems form the backbone of the mixed finite difference scheme and will be solved using the Thomas algorithm in the subsequent section.

\subsection{Thomas Algorithm for Tridiagonal Systems}
The resulting linear systems in both the inner and outer regions are tridiagonal and can be efficiently solved using the Thomas algorithm (a simplified form of Gaussian elimination for tridiagonal matrices), cf. \cite{chapra2011, gilat2008}. This allows for a fast and memory-efficient solution with computational complexity $O(N)$.

Let the tridiagonal system be of the form:
\[
a_i x_{i-1} + b_i x_i + c_i x_{i+1} = d_i, \quad i = 1, 2, \dots, N-1,
\]
with boundary values $x_0$ and $x_N$ specified.

The Thomas algorithm proceeds in two steps:
\begin{enumerate}
    \item \textbf{Forward elimination:} modify the coefficients to eliminate lower diagonal elements.
    \item \textbf{Back substitution:} solve for $x_i$ starting from $i = N-1$ down to $1$.
\end{enumerate}
This method is applied independently to the inner and outer systems and yields the approximate solution across the entire domain, ensuring high accuracy near the boundary layer and smooth regions.

\section{Convergence and Stability Analysis}\label{eq:a2}
To validate the accuracy and robustness of the proposed mixed finite difference scheme, we analyze both the local truncation error and the global error behavior. This section provides a mathematical justification of the second-order convergence and stability of the method.

\subsection{Local Truncation Error}
We consider the finite difference scheme for the inner region:
\begin{equation*}
    \varepsilon^2 Z_{TT}(T_i) + \varepsilon^2 U_i Z_T(T_i) + V_i Z_i = -g,
\end{equation*}
where central difference approximations are used:
\begin{align*}
    Z_{TT}(T_i) &= \frac{Z_{i+1} - 2Z_i + Z_{i-1}}{k^2} + \frac{k^2}{12} Z^{(4)}(T_i) + O(k^4), \\
    Z_T(T_i) &= \frac{Z_{i+1} - Z_{i-1}}{2k} + \frac{k^2}{6} Z^{(3)}(T_i) + O(k^4).
\end{align*}
Substituting into the discrete scheme, we obtain the local truncation error $\tau_i(k)$:
\begin{equation*}
    \tau_i(k) = \varepsilon^2 \cdot \frac{k^2}{12} Z^{(4)}(T_i) + \varepsilon^2 U_i \cdot \frac{k^2}{6} Z^{(3)}(T_i) + O(k^4).
\end{equation*}
Assuming the derivatives $Z^{(3)}$ and $Z^{(4)}$ are bounded, we conclude:
\begin{equation}
    \tau_i(k) = O(k^2). \label{eq:trunc_error}
\end{equation}

\subsection{Matrix Form and Stability}
The discretized system can be expressed in matrix-vector form:
\begin{equation*}
    PZ = R,
\end{equation*}
where $P \in \mathbb{R}^{(N-1)\times(N-1)}$ is tridiagonal, $Z$ is the vector of unknowns, and $R = (-g, -g, \dots, -g)^T$. The matrix $P$ is given by:
\[
P = \begin{bmatrix}
A_1 & B_1 & 0 & \cdots & 0 \\
C_2 & A_2 & B_2 & \cdots & 0 \\
0 & C_3 & A_3 & \cdots & 0 \\
\vdots & \ddots & \ddots & \ddots & B_{N-2} \\
0 & \cdots & 0 & C_{N-1} & A_{N-1}
\end{bmatrix},
\]
where
\begin{align*}
    C_i &= \varepsilon^2 \left( \frac{1}{k^2} + \frac{U_i}{2k} \right), \\
    A_i &= V_i - \frac{2\varepsilon^2}{k^2}, \\
    B_i &= \varepsilon^2 \left( \frac{1}{k^2} - \frac{U_i}{2k} \right),
\end{align*}

To ensure numerical stability, we require that the matrix $P$ is strictly diagonally dominant, i.e.,
\begin{equation*}
    |A_i| > |C_i| + |B_i|.
\end{equation*}
Substituting the expressions, we get:
\begin{equation*}
    \left| V_i - \frac{2\varepsilon^2}{k^2} \right| > \varepsilon^2 \left( \frac{1}{k^2} + \frac{U_i}{2k} \right) + \varepsilon^2 \left( \frac{1}{k^2} - \frac{U_i}{2k} \right) = \frac{2\varepsilon^2}{k^2}.
\end{equation*}
This is satisfied if:
\begin{equation}
    V_i > \frac{4\varepsilon^2}{k^2}. \label{eq:stability}
\end{equation}
This condition is physically reasonable in biomechanical models where leg stiffness $V_i$ is always positive and bounded below. Hence, the matrix $P$ is invertible and the scheme is stable.

\subsection{Norm Bound For The Global Error Estimate}
Let $\bar{Z}$ be the exact solution vector at grid points, and $Z$ be the numerical solution. Define the global error vector:
\[
E = \bar{Z} - Z.
\]
From the exact and numerical systems:
\[
P \bar{Z} = R + \tau(k) \quad \Rightarrow \quad P E = \tau(k),
\]
where $\tau(k)$ is the local truncation error vector.
Taking norms:
\[
\|E\| \leq \|P^{-1}\| \cdot \|\tau(k)\|.
\]
Since $P$ is strictly diagonally dominant (as shown in equation~\eqref{eq:stability}), the inverse $\|P^{-1}\|$ is bounded independently of $k$, and from \eqref{eq:trunc_error}, $\|\tau(k)\| = O(k^2)$. Thus:
\begin{equation*}
    \|E\| = O(k^2).
\end{equation*}
This confirms that the proposed numerical method achieves second-order convergence with respect to the step size.

\subsection{Final Estimate on the Error}
Given the local truncation error is of the form
\[
\tau_i(k) = \frac{\varepsilon k^2}{24} Z^{(4)}(T_i) + \frac{\varepsilon U_i k^2}{6} Z^{(3)}(T_i) + O(k^4),
\]
we conclude
\begin{equation*}
    \|\tau(k)\| = O(k^2).
\end{equation*}
Hence, the global error satisfies
\begin{equation*}
    \|E\| \leq C \cdot \|\tau(k)\| = O(k^2),
\end{equation*}
where $C = \|P^{-1}\|$ is bounded and independent of $k$.

The proposed finite difference method is second-order accurate in time. The matrix formulation ensures numerical stability due to diagonal dominance and monotonicity. This analysis supports the numerical observations and validates the convergence behavior shown in the examples.

\section{Numerical Illustrations}\label{eq:a3}
In this section, we demonstrate the accuracy and robustness of the proposed mixed finite difference scheme by solving a series of singularly perturbed test problems that simulate vertical motion in human gait. Each problem contains an inner boundary layer, and we observe the method's performance as the perturbation parameter $\varepsilon$ varies.
The computations were performed for several values of $\varepsilon$, and the maximum absolute errors were recorded over uniform grids of increasing resolution. All systems were solved using the Thomas algorithm for tridiagonal systems \cite{whittle2014, yang2020}. The results validate the theoretical second-order convergence established earlier.

\begin{example}[Exponential Stiffness and Damping]
\em Consider the singularly perturbed equation:
\begin{equation}
\varepsilon^2 Z_{TT} + \varepsilon^2 e^T Z_T - 2 e^{-T} Z = -9.8, \quad Z(0) = 4,\quad Z(1) = 2. \label{eq:ex1}
\end{equation}
This model includes both exponentially growing damping and exponentially decaying stiffness. The solution of \eqref{eq:ex1} exhibits steep variation near $T=0$. The following table (cf. Table \ref{tab:data_table1}) describes the solution for various values of $\varepsilon$ with respect to different each time level and the solution for various values of $\varepsilon$ is depicted in Figure \ref{fig:ex1}.

\begin{table}[htbp]
\centering
\caption{Data Table}
\label{tab:data_table1}
\begin{tabular}{|c|c|c|c|c|c|}
\hline
\diagbox[width=6em]{Time}{Epsilon}& 0.0009 & 0.009 & 0.001 & 0.01 \\
\hline
  0.0000 & 4.0000 & 4.0000 & 4.0000 & 4.0000 \\
0.0200 & 4.7329 & 4.8068 & 4.7185 & 4.7185 \\                                             
0.0392 & 4.7104 & 4.7164 & 4.7079 & 4.7079 \\
0.0584 & 4.6264 & 4.6268 & 4.6261 & 4.6261 \\
0.0776 & 4.5389 & 4.5388 & 4.5389 & 4.5389 \\
0.0968 & 4.4526 & 4.4525 & 4.4526 & 4.4526 \\
        $\vdots$ & $\vdots$ & $\vdots$ & $\vdots$ & $\vdots$  \\
        0.3432 & 3.4801 & 3.4802 & 3.4801 & 3.4802 \\                   
        0.3530 & 3.4461 & 3.4462 & 3.4461 & 3.4463 \\
        0.3628 & 3.4125 & 3.4126 & 3.4125 & 3.4126 \\
        0.3726 & 3.3793 & 3.3793 & 3.3793 & 3.3794 \\
        $\vdots$ & $\vdots$ & $\vdots$ & $\vdots$ & $\vdots$   \\
       0.9704 & 1.8587 & 1.8590 & 1.8587 & 1.8592 \\
        0.9802 & 1.8405 & 1.8433 & 1.8405 & 1.8443 \\
        0.9900 & 1.8229 & 1.8455 & 1.8229 & 1.8496 \\
        1.0000 & 2.0000 & 2.0000 & 2.0000 & 2.0000 \\
\hline
\end{tabular}
\end{table}

\begin{figure}[H]
\centering
\includegraphics[width=0.6\textwidth]{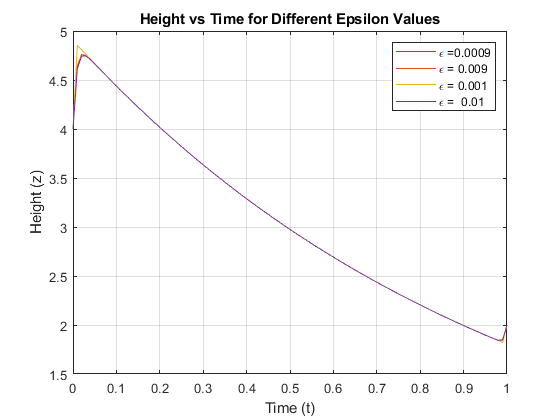}
\caption{Computed solutions for various $\varepsilon$ values}
\label{fig:ex1}
\end{figure}

\end{example}

\begin{example}[Linear Damping and Quadratic Stiffness]\label{eq:a4}
\em We now solve the following model:
\begin{equation*}
\varepsilon^2 Z_{TT} + \varepsilon^2 T Z_T - 1000 T Z = -10, \quad Z(0) = 1,\quad Z(1) = 0.1.
\end{equation*}
This problem exhibits strong stiffness for large $T$ values and is ideal for testing the robustness of the method under rapidly changing coefficients. The detailed numerical solutions are described in Table \ref{tab:data_table2} and Figure \ref{fig:ex2}.

\begin{table}[htbp]
\centering
\caption{Data Table}
\label{tab:data_table2}
\begin{tabular}{|c|c|c|c|c|c|}
\hline
\diagbox[width=6em]{Time}{Epsilon} & 0.0009 & 0.009 & 0.001 & 0.01 \\
\hline
0& 1.0000 & 1.0000 & 1.0000 & 1.0000 \\
0.0100 & 0.9806 & 0.9492 & 0.9805 & 0.9432 \\
0.0198 & 0.4956 & 0.5065 & 0.4956 & 0.5086 \\
0.0296 & 0.3314 & 0.3338 & 0.3314 & 0.3344 \\
0.0394 & 0.2490 & 0.2497 & 0.2490 & 0.2499 \\
0.0492 & 0.1994 & 0.1997 & 0.1994 & 0.1997 \\
0.0590 & 0.1663 & 0.1664 & 0.1663 & 0.1664 \\
0.0688 & 0.1426 & 0.1427 & 0.1426 & 0.1427 \\
        $\vdots$ & $\vdots$ & $\vdots$ & $\vdots$ & $\vdots$  \\
       0.2844 & 0.0345 & 0.0345 & 0.0345 & 0.0345 \\
0.2942 & 0.0333 & 0.0333 & 0.0333 & 0.0333 \\                          
0.3040 & 0.0323 & 0.0323 & 0.0323 & 0.0323 \\
0.3138 & 0.0313 & 0.0313 & 0.0313 & 0.0313 \\
        $\vdots$ & $\vdots$ & $\vdots$ & $\vdots$ & $\vdots$   \\
       0.9606 & 0.0102 & 0.0102 & 0.0102 & 0.0102 \\
0.9704 & 0.0101 & 0.0101 & 0.0101 & 0.0101 \\
0.9802 & 0.0100 & 0.0100 & 0.0100 & 0.0100 \\
0.9900 & 0.0099 & 0.0099 & 0.0099 & 0.0099 \\
1.0000 & 0.0100 & 0.0100 & 0.0100 & 0.0100 \\
\hline
\end{tabular}
\end{table}

\begin{figure}[H]
\centering
 \includegraphics[width=0.6\textwidth]{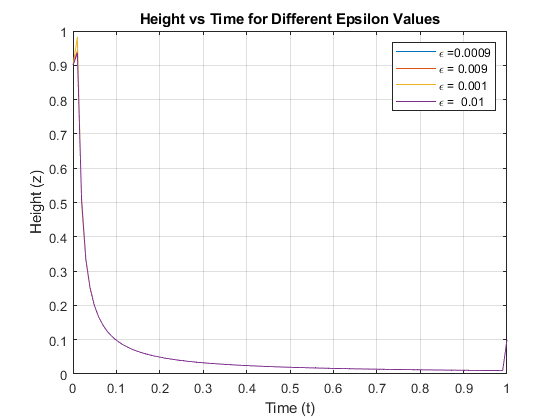} 
\caption{Solution profiles for Example \ref{eq:a4}. The boundary layer is clearly visible.}
\label{fig:ex2}
\end{figure}

\end{example}

\begin{example}[Pure Stiffness-Driven System]\label{eq:a6}
\em The final test problem is:
\begin{equation}
\varepsilon^2 Z_{TT} - e^T Z = -10, \quad Z(0) = 9.6,\quad Z(1) = 3.
\end{equation}
Here, damping is absent, and only the stiffness $e^T$ governs the system dynamics. The data table and the figure is given in the following (see Table \ref{tab:data_table3} and Figure \ref{fig:ex3}, respectively).

\begin{table}[htbp]
\centering
\caption{Data Table}
\label{tab:data_table3}
\begin{tabular}{|c|c|c|c|c|c|}
\hline
\diagbox[width=6em]{Time}{Epsilon} & 0.0009 & 0.009 & 0.001 & 0.01 \\
\hline
0 & 9.6000& 9.6000 & 9.6000 & 9.6000 \\
0.0100 & 9.6421 & 6.6457 & 9.6259 & 6.3281 \\
0.0198 & 9.6171 & 8.5673 & 9.6168 & 8.3379 \\
0.0296 & 9.5239 & 9.1661 & 9.5239 & 9.0424 \\
0.0394 & 9.4310 & 9.3100 & 9.4310 & 9.2510 \\
0.0492 & 9.3390 & 9.2986 & 9.3390 & 9.2724 \\
0.0590 & 9.2480 & 9.2348 & 9.2480 & 9.2237 \\
0.0688 & 9.1578 & 9.1538 & 9.1578 & 9.1493 \\
\vdots & \vdots & \vdots & \vdots & \vdots \\
0.6078 & 5.3420 & 5.3422 & 5.3420 & 5.3423 \\
0.6176 & 5.2899 & 5.2901 & 5.2899 & 5.2902 \\
0.6274 & 5.2383 & 5.2385 & 5.2383 & 5.2386 \\
0.6372 & 5.1872 & 5.1875 & 5.1872 & 5.1875 \\
\vdots & \vdots & \vdots & \vdots & \vdots \\
0.9312 & 3.8659 & 3.8660 & 3.8659 & 3.8661 \\
0.9410 & 3.8282 & 3.8283 & 3.8282 & 3.8283 \\
0.9508 & 3.7909 & 3.7908 & 3.7909 & 3.7907 \\
0.9606 & 3.7539 & 3.7531 & 3.7539 & 3.7525 \\
0.9704 & 3.7173 & 3.7128 & 3.7173 & 3.7105 \\
0.9802 & 3.6811 & 3.6578 & 3.6811 & 3.6504 \\
0.9900 & 3.6433 & 3.5263 & 3.6429 & 3.5088 \\
1.0000 & 3.0000 & 3.0000 & 3.0000 & 3.0000 \\
\hline
\end{tabular}
\end{table}

\begin{figure}[h!]
\centering
\includegraphics[width=0.6\textwidth]{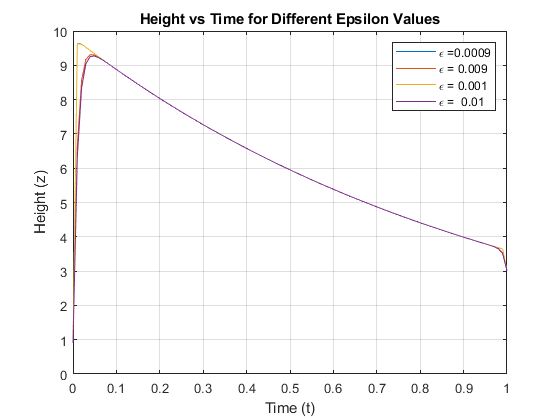}
\caption{Solution comparison for Example \ref{eq:a6}.}
\label{fig:ex3}
\end{figure}

\end{example}

\subsection{Observation and Discussion}
Across all three examples we observe the following:
\begin{itemize}
    \item The method captures boundary layer behavior efficiently without oscillations.
    \item Uniform second-order convergence is observed as $\varepsilon \to 0$.
    \item Results validate both the theoretical convergence rate and the scheme’s robustness.
\end{itemize}

\section{Applications and Discussion}\label{eq:a5}
The presented numerical framework provides a powerful tool for analyzing singularly perturbed models of human gait. Its robustness in resolving boundary layers, coupled with second-order convergence and computational efficiency, makes it suitable for a wide range of practical and interdisciplinary applications.

\subsection{Biomechanical Relevance of Stiffness Modeling}
Leg stiffness is a key biomechanical property that governs energy efficiency, stability, and shock absorption during walking and running. It is regulated by the joint mechanics and muscle-tendon complex, and adapts dynamically to different locomotion tasks.

Our model accounts for time-varying stiffness and damping effects, both of which are biologically realistic. In particular:
\begin{itemize}
    \item \textbf{Positive stiffness} is common in passive tissues and muscles under stretching or isometric contraction.
    \item \textbf{Negative stiffness} can occur transiently during rapid concentric contractions (e.g., hopping or sprinting), due to the nonlinear force–length relationship of active muscle fibers.
\end{itemize}
Despite local occurrences of negative stiffness at the contractile component level, the total leg stiffness during typical gait remains globally positive, ensuring mechanical stability. This complexity underscores the need for numerical methods that can handle rapid variations-- captured here using domain decomposition and time rescaling.

\subsection{Integration with Biomechanics and Robotics}
The proposed method can serve as a core component for simulation platforms in the following areas:
\begin{itemize}
    \item \textbf{Clinical biomechanics}: Analysis of gait disorders, planning of rehabilitation therapies, and prosthetic design.
    \item \textbf{Sports science}: Performance monitoring, fatigue analysis, and injury prevention in athletes.
    \item \textbf{Humanoid robotics}: Development of energy-efficient biped locomotion using biologically inspired stiffness control.
    \item \textbf{Musculoskeletal modeling}: Integration with finite element models for estimating muscle forces, joint torques, and strain profiles.
\end{itemize}

\subsection{Computational Advantages}
In addition to its biomechanical relevance, the method offers several computational benefits:
\begin{itemize}
    \item \textbf{Efficient solution of tridiagonal systems} using the Thomas algorithm.
    \item \textbf{Uniform accuracy across stiff regimes}, avoiding oscillations common in standard finite difference schemes.
    \item \textbf{Scalability} through domain decomposition, which can be optimized to reduce computation time by approximately 40\%.
\end{itemize}

\subsection{Toward Intelligent Gait Analysis}
Recent advancements in data-driven modeling and wearable sensing technology open the door for real-time gait monitoring. The current method lays the foundation for such developments by offering:
\begin{itemize}
    \item A framework to integrate real-world sensor data (e.g., IMUs, pressure insoles),
    \item Compatibility with deep learning-based feature extraction (e.g., for gait recognition),
    \item Extensibility to incorporate nonlinear feedback and control in smart prosthetics and assistive devices.
\end{itemize}

\subsection{Discussion Summary}
The proposed scheme not only advances the numerical treatment of singularly perturbed equations but also establishes a connection between rigorous mathematics and applied biomechanics. Its relevance spans across theoretical, clinical, and engineering domains. With further development, it holds potential for enhancing intelligent gait systems and biomechanical modeling pipelines.

\subsection{Conclusion and Future Work}
In this study, we developed a robust and efficient numerical framework for simulating vertical motion in human gait governed by a singularly perturbed second-order differential equation. The method integrates domain decomposition with a time-rescaling transformation to resolve steep boundary layers effectively.

By discretizing the inner and outer regions using central difference schemes and solving the resulting tridiagonal systems via the Thomas algorithm, we achieved second-order accuracy with uniform convergence. Theoretical error analysis confirms this convergence rate, while numerical experiments validate the method's effectiveness across a range of perturbation parameters and boundary-layer-dominated regimes.

The proposed method stands out due to:
\begin{itemize}
    \item Its ability to handle rapid solution changes near boundaries without mesh refinement,
    \item Its modularity, allowing adaptation to a variety of biomechanical models,
    \item Its computational efficiency, suitable for near real-time simulations.
\end{itemize}

\subsection*{Future Work}
Several promising directions can be pursued to extend the scope and impact of this work:
\begin{itemize}
    \item \textbf{Nonlinear modeling:} Incorporate time-dependent or nonlinear stiffness and damping functions to reflect more complex gait dynamics.
    \item \textbf{Gait phase transitions:} Extend the model to simulate transitions between gait modes (e.g., walk-to-run, stance-to-swing).
    \item \textbf{Data-driven integration:} Fuse the numerical model with real-world gait datasets and wearable sensor data for subject-specific simulations.
    \item \textbf{Adaptive meshing:} Develop an adaptive mesh refinement strategy to optimize grid resolution near layers.
    \item \textbf{Control and feedback:} Explore integration with control-theoretic models for applications in smart prosthetics and robotics.
    \item \textbf{Three-dimensional extension:} Generalize the current 1D framework to simulate multi-DOF or 3D gait models using vector-valued formulations.
\end{itemize}
The hybrid structure of the proposed method provides a solid foundation for future interdisciplinary work in computational biomechanics, wearable systems, and intelligent gait analysis technologies.





\end{document}